\documentclass{amsart}
\newtheorem{Theorem}{Theorem}[section]
\newtheorem{Definition}[Theorem]{Definition}
\newtheorem{Proposition}[Theorem]{Proposition}

\newtheorem{Lemma}[Theorem]{Lemma}
\newtheorem{Corollary}[Theorem]{Corollary}
\theoremstyle{remark}

\newtheorem{Example}[Theorem]{Example}

\def\il{\int\limits_}

\def\eps{\epsilon}

\def\ovr{\overline}

\def\om{\omega}
\def\Om{\Omega}

\def\gm{\gamma}
\def\th{\theta}

\def\dl{\delta}

\def\bd{\partial}
\def\lm{\lambda}

\def\si{\sigma}

\def\sm{\setminus}
\def\sbs{\subset}

\def\cl{\operatorname{cl}}

\def\supp{\operatorname{supp}}

\def\ess{\operatorname{ess}}

\def\be{\begin{enumerate}}
\def\ee{\end{enumerate}}
\def\bT{\begin{Theorem}}
\def\eT{\end{Theorem}}
\def\bP{\begin{Proposition}}
\def\eP{\end{Proposition}}
\def\bD{\begin{Definition}}
\def\eD{\end{Definition}}
\def\bE{\begin{Example}}
\def\eE{\end{Example}}
\def\bL{\begin{Lemma}}
\def\eL{\end{Lemma}}
\def\bC{\begin{Corollary}}
\def\eC{\end{Corollary}}

\def\A{{\mathcal A}}

\def\J{{\mathcal J}}

\def\M{{\mathcal M}}

\def\rS{{\mathcal S}}
\def\T{{\mathcal T}}

\begin{document}
\title{Dirichlet problems for plurisubharmonic functions on compact sets}
\author{Evgeny A. Poletsky and Ragnar Sigurdsson}
\begin{abstract} In this paper we solve the Dirichlet problems for
different classes of plurisubharmonic functions on compact sets in
$\mathbb C^n$ including continuous, pluriharmonic and maximal
functions.
\end{abstract}
\thanks{The first author was supported by the NSF Grant DMS-0900877.}
\address{ Department of Mathematics,  Syracuse University, \newline
215 Carnegie Hall, Syracuse, NY 13244} \email{eapolets@syr.edu}
\address{ Science Institute, University of Iceland,
\newline Dunhaga 3, IS-107 Reykjavik,
Iceland} \email{ragnar@hi.is} \maketitle
\section{Introduction}
\par To define plurisubharmonic functions on a compact set $X$ one,
first, introduces the sets of Jensen measures $\J_z(X)$ on $X$
with barycenter $z\in X$ as the weak-$*$ limits of Jensen
measures on neighborhoods of $X$. Then plurisubharmonic
functions are defined as functions with appropriate topological
properties and satisfying subaveraging inequality with respect
to these measures.
\par Compact sets carrying plurisubharmonic functions appear
rather naturally in pluripotential theory and the theory of
uniform algebras but not too much is known about their
structure. In particular, it is not known how flexible the
space of plurisubharmonic functions on compact sets is.
\par In this paper we approach this problem through a study of
Dirichlet problems for different classes of plurisubharmonic
functions on compact sets in $\mathbb C^n$. To pose a Dirichlet
problem we need a notion of a boundary. Our claim is that the
natural boundary for plurisubharmonic functions is the closure
$B_X$ of the set $O_X=\{z\in X:\,\J_z(X)=\{\dl_z\}\}$.
\par The fact that it is natural is confirmed in Section
\ref{S:lsp} where we prove a generalization of the Littlewood
subordination principle: Any Jensen measure is subordinated to
a Jensen measure supported by $B_X$. Another confirmation comes
in Section \ref{S:orcs} where we show that if $O_X=B_X$ then
any continuous function on $B_X$ can be extended to $X$ as a
continuous plurisubharmonic function.
\par In Section \ref{S:ps} we describe compact sets on which any
continuous function on $B_X$ can be extended as a pluriharmonic
function. We prove that such compact sets are rather rare: The
Dirichlet problem for pluriharmonic functions can be solved if
and only if for every $z\in X$ there is only one Jensen measure
in $\J_z(X)$ supported by $B_X$.
\par In Section \ref{S:msdp} we look for maximal solutions of the
Dirichlet problem. An example in this section demonstrates
that, in general, this problem does not admit a continuous
solution. However, we prove that it always admits a lower
semicontinuous solution with subaveraging property. We also
give in this section a description of maximal functions in
terms of Jensen measures.
\par In Section \ref{S:homf} we address the problem of covering
the domain of a maximal function $u$ with compact sets such
that the restriction of $u$ to each of these sets is harmonic.
A result of Bedford and Kalka \cite{BK} established the
existence of such a covering by complex manifolds in a smooth,
non-degenerate case. However, there are examples of continuous
maximal functions (see, for example, \cite{D} and \cite{DL})
where there are no coverings with analytic structure. We give a
definition of harmonic functions on compact sets in terms of
Jensen measures and provide sufficient conditions for the
existence of coverings above.
\par This paper was written when the first author was visiting the
Department of Mathematics of the University of Iceland and the
second author was visiting the Department of Mathematics of
Syracuse University. We express our gratitude to both
departments for their hospitality.
\section{Basic notions and facts}\label{S:bnf}
\par Let ${\mathbb D}$ be the unit disk in ${\mathbb C}$ and
${\mathbb T}=\bd{\mathbb D}$. For a set $E\sbs\ovr{\mathbb D}$ we
consider the infimum $v(z)$ of all positive superharmonic
functions on ${\mathbb D}$, extending continuously on
$\ovr{\mathbb D}$ and greater than 1 on $E$. Then the function
\[\om(z,E,{\mathbb D})=\liminf\limits_{w\to z}v(z)\] is called the
{\it harmonic measure} of the set $E$ with respect to ${\mathbb
D}$. The function $\om(z,E,{\mathbb D})$ is superharmonic on
${\mathbb D}$ and harmonic on ${\mathbb D}$ when $E\sbs{\mathbb
T}$.
\par Let us denote by $\lambda$ the normalized arc length measure
on ${\mathbb T}$. Let $f$ be a bounded holomorphic mapping of the
unit disk ${\mathbb D}$ into ${\mathbb C}^n$. Since $f$ has radial
limit values almost everywhere on ${\mathbb T}$, we may consider
$f$ as a Borel measurable mapping of the closed unit disk
$\ovr{\mathbb D}$. The push-forward
\[\mu_f(E)=\lm(f^{-1}(E)\cap{\mathbb T})=
\om(0,f^{-1}(E)\cap{\mathbb T},{\mathbb D})\] of $\lm$ by the
restriction of $f$ to ${\mathbb T}$ is a regular Borel measure on
${\mathbb C}^n$. The point $z_f=f(0)$ is uniquely determined as
the point such that
\[z_f=\int z\,\mu_f(dz).\]
\par Let $L=\{f_j\}$ be a sequence of uniformly bounded holomorphic
mappings of the unit disk ${\mathbb D}$ into ${\mathbb C}^n$.
{\it The cluster limit set} $\cl L$, or, short, {\it the
cluster } of the sequence $L$ is the set of all points
$z\in{\mathbb C}^n$ such that for every $r>0$ and infinitely
many $j$ the sets $f_j({\mathbb D})\cap B(z,r)$ are nonempty,
where $B(z,r)$ is the open ball of radius $r$ centered at $z$
in ${\mathbb C}^n$. The sequence $L$ is said to be {\it
weak-$*$ converging } if the measures $\mu_{f_j}$ converge
weak-$*$ to a measure $\mu_L$, i.e.,
\[\lim\limits_{j\to\infty}\il {{\mathbb C}^n}\phi(z)\,\mu_{f_j}(dz)=
\il {{\mathbb C}^n}\phi(z)\,\mu_L(dz)\]for every continuous
function $\phi$ on ${\mathbb C}^n$. By the Alouglou theorem
\cite[V.4.2]{C}, every sequence of uniformly bounded holomorphic
mappings contains a weak-$*$ converging subsequence. If
$L=\{f_j\}$ is weak-$*$ converging then the sequence $\{z_{f_j}\}$
converges to a point $z_L$, which we will call {\it the center }
of $L$.
\par A point $z\in\cl L$ is called {\it totally essential} if
\[\Om(z_L,V,L)=\liminf\limits_{j\to\infty}\om(0,f_j^{-1}(V),{\mathbb D})>0\]
for every open set $V$ containing $z$. Clearly, the set $\ess L$
of totally essential points is closed. Other points in $\cl L$ are
called {\it nonessential.}
\par If a sequence  $L=\{f_j\}$ is weak-$*$ converging, then
$\supp L=\supp \mu_L$ is compact, every point $z\in\supp L$ is
totally essential, and $\supp L=\{z_L\}$ if and only if $\ess
L=\{z_L\}$. Such a sequence is called  {\it totally perfect }
if all points of its cluster are totally essential. All
subsequences of a totally perfect sequence have the same
cluster and are totally perfect.
\par By \cite[Th. 2.2]{P2} if $L=\{f_j\}$ is a weak-$*$ converging
sequence then for every $z_0\in\cl L$ there is a totally
perfect sequence $M$ such that $z_M=z_L$, $z_0\in\cl M$,
$\mu_M=\mu_L$, $\cl M\sbs\cl L$, and $\cl M$ contains all
totally essential points of $L$.
\par Let $X$ be a compact set in ${\mathbb C}^n$. We denote by $\T(X)$
the set of all totally perfect sequences $L$ such that $\cl L\sbs
X$ and by $\M(X)$ the set of all measures $\mu_L$, $L\in\T(X)$. If
a point $z\in X$ we denote by $\T_z(X)$ the set of all $L\in\T(X)$
such that $z_L=z$. Let $\M_z(X)$ be the set of measures $\mu_L$,
$L\in\T_z(X)$.
\par A function $u$ on $X$ is {\it plurisubharmonic } if it is
upper semicontinuous and $u(z)\le\mu_L(u)$ for all $L\in\T_z(X)$.
We allow plurisubharmonic functions to assume $-\infty$ as its
value and we say that a function on $X$ is continuous if it
continuously maps $X$ into $[-\infty,\infty)$. We will denote the
space of plurisubharmonic functions on $X$ by $PSH(X)$. Let
$PSH^c(X)=PSH(X)\cap C(X)$. A function $u(z)$ on $X$ is {\it
plurisuperharmonic } if $-u\sbs PSH(X)$.
\par It was proved in \cite[Cor. 3.1]{P2} that every function
$u\in PSH^c(X)$ is the limit of an increasing sequence of
continuous plurisubharmonic functions defined in neighborhoods of
$X$.
\par {\it A Jensen measure on $X$ with barycenter $z_0\in X$ } is a
regular nonnegative Borel measure $\mu$ supported by $X$ such that
$\mu(X)=1$ and $u(z_0)\le\mu(u)$ for every plurisubharmonic
function $u$ on $X$. We denote the set of such measures by
$\J_{z_0}(X)$. Clearly, this set is convex and weak-$*$ compact.
By \cite[Thm. 3.2]{P2} $\M_{z_0}(X)=\J_{z_0}(X)$.
\par Jensen measures enjoy the following compactness property.
\bL\label{L:lpjm} Suppose that $\{X_j\}$ is a decreasing sequence of compact
sets in a closed ball $B\sbs\mathbb{C}^n$, $X=\cap X_j$ and a
point $z_j\in X_j$. Let $\{\mu_j\}$ be a sequence of measures such
that $\mu_j\in\J_{z_j}(X_j)$. Then any limit point $\mu$ of this
sequence in $C^*(B)$ belongs to $\J_{z_0}(X)$, where $z_0$ is a
limit point of $\{z_j\}$.
\eL
\begin{proof} Let $L_j=\{f_{jk}\}\in\T_{z_0}(X_j)$ be totally perfect
sequences such that $\mu_{L_j}=\mu_j$. We may assume that the
measures $\mu_j$ converge weak-$*$ to $\mu$ and $f_{jk}(\mathbb
D)$ lies in the $1/j$-neighborhood of $X$. Since the space $C(B)$
is separable the weak-$*$ topology on the unit ball of $C^*(B)$ is
metrizable. Hence we can choose a sequence $L=\{f_{jk_j}\}$ such
that the measures $\mu_{f_{jk_j}}$ converge weak-$*$ to $\mu$ and
points $f_{jk_j}(0)$ converge to $z_0$. Clearly, $\cl L\sbs X$,
$z_L=z_0$ and, by \cite[Thm. 2.2]{P2}, we may assume that
$L\in\T_{z_0}(X)$. Thus $\mu=\mu_L\in\M_{z_0}(X)=\J_{z_0}(X)$.
\end{proof}
\par If $\phi$ is a function on $X$ then we denote by $E_X\phi$
the upper envelope of all continuous plurisubharmonic functions
$u\le\phi$ on $X$. The main tool in dealing with Jensen
measures is
\bT[Edwards' Theorem]\label{T:elscf} Let $\phi$ be a lower semicontinuous
function on a compact set $X\sbs\mathbb C^n$. Then
\[E_X\phi(z)=\inf\{\mu(\phi):\,\mu\in\J_z(X)\}.\]
Moreover, the infimum is attained.
\eT
\par Formally, the infimum above should be taken over Jensen
measures with respect to $PSH^c(X)$. But it was proved in
\cite{P2} that if $\mu$ is a regular nonnegative Borel measure
on $X$ and $u(z_0)\le\mu(u)$ for every continuous
plurisubharmonic function $u$ on $X$, then $\mu$ is a Jensen
measure.
\par For a point $z\in X$ we define the set $I_z$ as the set of all
points $w\in X$ such that $w\in\cl L$ for some $L\in\T_z(X)$.
Let $O_X$ be the set of all points $z\in X$ such that
$I_z=\{z\}$ or, what is the same, $\J_z(X)=\{\dl_z\}$. It was
shown in \cite{P2} that the set $O_X$ is nonempty. We denote by
$B_X$ the closure of $O_X$. We shall call $B_X$ {\it the
potential boundary } of $X$. We let $\J^b_z(X)$ be the set of
all $\mu\in\J_z(X)$ such that $\supp\mu\sbs B_X$. By \cite{P2}
$B_X$ is the minimal closed set such that for all points $z\in
X$ there is $L\in\T_z(X)$ with $\supp L\sbs B_X$. Thus
$\J^b_z(X)$ is non-empty for all $z\in X$.
\par The following simple lemma says a little bit more about
$B_{\cl L}$.
\bL\label{L:pbc} Let $X$ be a compact set in $\mathbb C^n$,
$z_0\in X$, $L\in\T_{z_0}(X)$ and $\supp L\sbs B_X$. If $Y=\cl
L$ then $B_Y\sbs B_X$.
\eL
\begin{proof} By \cite[Lemma 2.2]{P2} for every $z\in Y$ there is
$M\in\T_z(Y)$ such that $\supp \mu_M\sbs\supp L$. Since $\supp
L\sbs B_X$ we see that any point $z\in\cl L\sm B_X$ does not
belong to $B_Y$.
\end{proof}
\section{The Littlewood subordination principle}\label{S:lsp}
\par For $\mu,\nu\in \J_z(X)$ we say that $\nu$ is {\it subbordinated }
to $\mu$ and denote it $\mu\succeq\nu$ if $\mu(u)\ge\nu(u)$ for
all $u\in PSH^c(X)$. Clearly, if $\mu\succeq\nu$ and
$\nu\succeq\eta$, then $\mu\succeq\eta$.
\bL\label{L:po} If $\mu,\nu\in\J_z(X)$, $\mu\succeq\nu$ and
$\nu\succeq\mu$, then $\mu=\nu$.
\eL
\begin{proof} Clearly, $\mu(u)=\nu(u)$ for every $u\in PSH^c(X)$. If a
function $\phi\in C^2(X)$, i.e. $\phi$ is a $C^2$ function
defined on a neighborhood of $X$, then $u(z)=\phi(z)+k|z|^2\in
PSH^c(X)$ when $k$ is sufficiently large. Since
$\phi(z)=u(z)-k|z|^2$ we see that $\mu(\phi)=\nu(\phi)$.
Approximating a function $\phi\in C(X)$ by functions from
$C^2(X)$ we derive that $\mu(\phi)=\nu(\phi)$. Hence $\mu=\nu$.
\end{proof}
\par It follows that $\mu\succeq\nu$ is a partial order on
$\J_z(X)$. A measure $\mu\in \J_z(X)$ is called {\it maximal}
if there is no $\nu\in \J_z(X)$ such that $\nu\succeq\mu$ and
$\nu\ne\mu$. It is known (see \cite[Lemma 4.1]{Ph} that for any
$\nu\in \J_z(X)$ there is a maximal $\mu\in \J_z(X)$ such that
$\mu\succeq\nu$.
\par  The theorem below is a version of the classical Littlewood subordination
principle. Since this principle has many different interpretations
we state and prove what we have in mind. Before we do it let us
recall that if a subharmonic function $\phi\not\equiv-\infty$ on
${\mathbb D}$ has a harmonic majorant, then (see \cite[Ex. II.17,
19]{Ga}) it has the least harmonic majorant $h$ and $\phi=h+u$,
where $u$ is a subharmonic function on $\mathbb D$ whose radial
limits are equal to $0$ a.e. on $\mathbb T$ and
\[\lim_{r\to1^-}\il0^{2\pi}u(re^{i\th})\,d\th=0.\]
If
\[\|h\|^p_p=\limsup_{r\to1^-}\il0^{2\pi}|h(re^{i\th})|^p\,d\th<\infty\] for some
$p>1$, then (see \cite[Thm. IX.2.3]{Go}) $h$ has radial limits
$h^*(e^{i\th})$ a.e. on $\mathbb T$,
\[h(z)=\frac1{2\pi}\il0^{2\pi}P(z,\th)h^*(e^{i\th})\,d\th,\]
where $P$ is the Poisson kernel, and
\[\|h\|^p_p=\il0^{2\pi}|h^*(re^{i\th})|^p\,d\th.\]
Consequently, $\phi$ has radial limits
$\phi^*(e^{i\th})=h^*(e^{i\th})$ a.e. on $\mathbb T$ and
\[\lim_{r\to1^-}\il0^{2\pi}\phi(re^{i\th})\,d\th=\il0^{2\pi}\phi^*(e^{i\th})\,d\th.\]
\bT\label{T:lmp} If a subharmonic function $\phi\not\equiv-\infty$ on
${\mathbb D}$ has a harmonic majorant $h$ with $\|h\|_p<\infty$
for some $p>1$,  then
\[\frac1{2\pi}\il0^{2\pi}(\phi\circ f)^*(e^{i\th})\,d\th\le
\frac1{2\pi}\il0^{2\pi}\phi^*(e^{i\th})\,d\th\]
for any holomorphic mapping $f:\,{\mathbb D}\to{\mathbb D}$ with
$f(z_0)=0$.
\eT
\begin{proof} First, we assume that $f$ is continuous up to the boundary and
$f(\ovr{\mathbb D})\sbs\mathbb D$. Then
\[\frac1{2\pi}\il0^{2\pi}\phi(f(e^{i\th}))\,d\th\le
\frac1{2\pi}\il0^{2\pi}h(f(e^{i\th}))\,d\th=h(0)=
\frac1{2\pi}\il0^{2\pi}\phi^*(e^{i\th})\,d\th.\] In the general
case, for $0<r<1$ we let $f_r(\zeta)=f(r\zeta)$, $|\zeta|\le1$.
Since $|h|^p$, $p\ge 1$, is a subharmonic function on $\mathbb D$
by the previous result
\[\il0^{2\pi}|h(f_r(e^{i\th}))|^p\,d\th\le\il0^{2\pi}|h^*(e^{i\th})|^p\,d\th.\]
Hence, $\|h\circ f\|_p<\infty$ and
\[\il0^{2\pi}(\phi\circ f)^*(e^{i\th})\,d\th=
\lim_{r\to1^-}\il0^{2\pi}\phi(f_r(e^{i\th}))\,d\th\le
\il0^{2\pi}\phi^*(e^{i\th})\,d\th.\]
\end{proof}
\par The following theorem is a generalization of this
principle.
\bT\label{T:lsp}  Let $X\subset{\mathbb C}^n$ be a compact set, $z_0\in X$
and $\mu\in \J_{z_0}(X)$.  Then there is a measure
$\nu\in\J_{z_0}^b(X)$ such that $\nu\succeq\mu$.
\eT
\begin{proof} If $f:\,{\mathbb D}\to{\mathbb C}^n$ is a holomorphic
mapping and $f({\mathbb D})\sbs B$, where $B$ is a closed ball in
${\mathbb C}^n$, then
\[\lim_{t\to1^-}\il0^{2\pi}\phi(f(te^{i\th}))\,d\th=
\il0^{2\pi}\phi(f(e^{i\th}))\,d\th\] for every function $\phi\in
C(B)$. This means that if $f_t(\zeta)=f(t\zeta)$, $0<t<1$, then
measures $\mu_{f_t}$ converge weak-$*$ to $\mu_f$ as $t\to1^-$.
\par Since the space $C(B)$ is separable, the unit ball in $C^*(B)$
is metrizable. Therefore, if $L=\{f_j\}$ is a weak-$*$ converging
sequence such that $f_j({\mathbb D})\sbs B$, then there is a
sequence $\{t_j\}$, $0<t_j<1$, such that the sequence $M$ of
$g_j(\zeta)=f(t_j\zeta)$ converges weak-$*$ to $\mu_L$. Moreover,
$\cl M\sbs\cl L$.
\par Let $L=\{f_j\}\in\T_{z_0}(X)$ be a weak-$*$ converging sequence
such that $\mu=\mu_L$. By the argument above we may assume that
the mappings $f_j$ are holomorphic in a neighborhood of
$\ovr{\mathbb D}$.
\par Let $V_j$ be a sequence of open sets such that $V_{j+1}\sbs\sbs V_j$
and $\cap V_j=X$. We take a decreasing sequence of open sets
$W_j\subset\subset V_j$, such that $B_X\subset W_j$ and $\cap
W_j=B_X$.
\par Take $\phi_j\in C(V_j)$ with $\supp \phi_j\subset \subset W_j$
and such that $\phi_j\geq -1$ and $\phi_j=-1$ on $B_X$.  Let
\[v_j(z)=E_{V_j}\phi_j(z)=
\sup\{v:\,v\text{ is plurisubharmonic  on } V_j\text{ and
}v\le\phi_j\}.\] It is known (see, \cite{P1}) that $v_j$ is
plurisubharmonic on $V_j$. By \cite[Cor. 4.2]{P2} for every $z\in
X$ there is $\sigma\in \J_z(X)$ such that $\supp \sigma\subset
B_X$. So $-1\le v_j(z)\le\si(\phi_j)=-1$. Hence $v_j=-1$ on $X$.
\par For any $j$ we let $U_j=\{z\in V_j :  v_j(z)<-1+\epsilon_j\}$,
where the sequence of positive numbers $\{\eps_j\}$ converges to
0. The sets $U_j$ are open subsets of $V_j$ containing $X$. Since
$\cl L\sbs X$, for every $j$ there exists $k_j$ such that
$f_{k}(\overline {\mathbb D})\subset U_j$ when $k\ge k_j$.
\par By Lemma 2.7 in \cite{LS} (see also equations (2.4) and (2.3)
there) for every $k\geq k_j$ there is a holomorphic mapping
$F_{jk}:\,\overline {\mathbb D}\times \overline{\mathbb D}\to U_j$
such that $F_{jk}(\zeta,0)=f_k(\zeta)$ for all $\zeta\in \mathbb
D$ and
\begin{equation}\begin{aligned}
&\frac1{4\pi^2}\int_0^{2\pi}\int_0^{2\pi}
\phi_j(F_{jk}(e^{i\th},e^{i\psi}))\,d\th d\psi\notag\\
&\le\frac1{2\pi}\int_0^{2\pi}v_j(f_k(e^{i\th}))\,d\th+
\eps_j\le-1+2\eps_j.\notag
\end{aligned}\end{equation}
\par We introduce the measures $\mu_{jk}$, $k\ge k_j$, by
\[
\mu_{jk}(\phi)= \frac1{4\pi^2}\int_0^{2\pi}\int_0^{2\pi}
\phi(F_{jk}(e^{i\th},e^{i\psi}))\,d\th d\psi
\] for all $\phi\in C_0(V_j)$.
Clearly $\mu_{jk}\in \J_{z_0}(V_j)$, where $z_0$ is the barycenter
of $\mu$.
\par If $u$ is plurisubharmonic on $V_j$, then
\begin{equation}\label{e:lsp1}\begin{aligned}
\mu_{jk}(u)&=\frac1{4\pi^2}\int_0^{2\pi}\int_0^{2\pi}
u(F_{jk}(e^{i\th},e^{i\psi}))\,d\th d\psi \geq
\frac1{2\pi}\int_0^{2\pi}
 u(F_{jk}(e^{i\th},0))\,d\th\\
&=\frac1{2\pi}\int_0^{2\pi}
 u(f_{k}(e^{i\th}))\,d\th=\mu_{f_k}(u).
\end{aligned}\end{equation}
Hence $\mu_{jk}(u)\geq \mu_{f_k}(u)$ for all plurisubharmonic
functions $u$ on $V_j$.
\par Since $\supp \phi_j\subset W_j$, we have
\begin{equation}\label{e:lsp2}
\mu_{jk}(W_j)\geq -\frac1{4\pi^2}\int_0^{2\pi}\int_0^{2\pi}
\phi_j(F_{jk}(e^{i\th},e^{i\psi}))\,d\th d\psi \geq 1-2\epsilon_j.
\end{equation}
\par By Alaoglou's theorem for each $j$ the sequence $\mu_{jk}$ has
a limit point $\mu_j$ which by Lemma \ref{L:lpjm} belongs to
$\J_{z_0}(\ovr V_j)$. By (\ref{e:lsp1}) $\mu_j(u)\geq
\mu_{f_k}(u)$ for every $k\ge k_j$ and every plurisubharmonic
function $u$ on $V_{j-1}$ and by (\ref{e:lsp2})
$\mu_j(W_{j-1})\ge1-2\epsilon_j$. If $u$ is a continuous
plurisubharmonic function on $V_{j-1}$, then
$\mu_j(u)\ge\lim_{k\to\infty}\mu_{f_k}(u)=\mu(u)$. Again by
Alaoglou's theorem and Lemma \ref{L:lpjm} the sequence $\{\mu_j\}$
has a limit point $\nu\in\J_{z_0}(X)$. Clearly, if $u$ is a
continuous plurisubharmonic function on a neighborhood of $X$,
then $\nu(u)\ge\mu(u)$. Moreover, since
$\mu_j(W_{j-1})\ge1-2\epsilon_j$, $\nu(B_X)=1$. Hence, 
$\supp\nu\sbs B_X$ and $\nu\in \J^b_{z_0}(X)$.
\par If $u$ is continuous and plurisubharmonic on $X$, then by \cite[Cor.
3.1]{P2}, there is an increasing sequence of continuous
plurisubharmonic functions $u_j$ each defined on some neighborhood
of $X$, such that $u=\lim u_j$ on $X$.  Since $\nu(u_j)\geq
\mu(u_j)$  for all $j$ we get $\nu(u)\geq \mu(u)$.
\end{proof}
\par Combining Lemma \ref{L:po} and Theorem \ref{T:lsp} we
obtain the following corollary.
\bC\label{C:smp} If $\mu$ is a maximal Jensen measure on $X$
then $\supp\mu\sbs B_X$.
\eC
\section{$O$-regular compact sets}\label{S:orcs}
\par We say that a point $z\in X$ is {\it a peak point }
if there is a plurisubharmonic function $u$ on $X$ such that
$u(z)=0$ and $u(w)<0$ when $w\ne z$. It was shown in \cite{P2}
that $O_X$ is the set of all peak points.
\bL\label{L:epp} If $z_0\in X$ is a peak point then there is a continuous
plurisubharmonic function $u$ on $X$ such that $u(z)=0$ and
$u(w)<0$ when $w\ne z$.
\eL
\par The proof is based on the following lemma
\bL\label{L:bc} If for every neighborhood $V$ of a point $z\in X$
there is a negative continuous plurisubharmonic function $u$ on
$X$ such that $u(z)=-1$ and $u(w)\le-2$ when $w\in X\sm V$, then
there is a continuous plurisubharmonic function $u$ on $X$ such
that $u(z)=0$ and $u(w)<0$ when $w\ne z$.
\eL
\par This lemma was proved in \cite{P2} (see Lemma 3.2) for plurisubharmonic
functions. But the proof holds without any changes for continuous
plurisubharmonic functions.
\par {\it Proof of Lemma \ref{L:epp}.} We just need to verify the
conditions of Lemma \ref{L:bc}. Let us take a continuous function
$\phi$ on $X$ such that $\phi(z_0)=-1$, $\phi(z)<-1$, $z\ne z_0$
and $\phi(z)\le-4$ on $X\sm V$, where $V$ is a neighborhood of
$z_0$. Let $v(z)=E_X\phi$. By \cite[Lemma 3.1]{P2} $v(z_0)=-1$,
$v<-1$ on $X\sm\{z_0\}$ and $v\le-4$ in $X\sm V$ and there is an
increasing sequence of continuous plurisubharmonic functions
$v_j$, defined on neighborhoods $V_j$ of $X$, converging pointwise
to $v$. Hence we can find a continuous plurisubharmonic function
$v_j$ on $X$ such that $-2\le v_j(z_0)\le -1$, $v_j<-1$ on $X$ and
$v_j\le-4$ on $X\sm V$. The function $u(z)=-v(z)/v(z_0)$ is
negative, $u(z_0)=-1$ and $u\le-2$ on $X\sm V$. Hence the
conditions of Lemma \ref{L:bc} hold. $\Box$
\par The classical maximum principle, stating that a non-constant plurisubharmonic
function cannot achieve maximum outside of the boundary, does not
hold for compact sets. For example, if
$X=\{(\zeta,t):\zeta\in\ovr{\mathbb D},t\in[-1,1]\}\sbs{\mathbb
C}^2$, then $PSH(X)$ consists of all upper semicontinuous
functions on $X$ subharmonic in $\zeta$ for each $t$ and
$B_X=\{(\zeta,t):\zeta\in\mathbb T,t\in[-1,1]\}$. So the function
$u(\zeta,t)=1-t^2$ is in $PSH(X)$ and attains its maximum at
$(0,0)$ away from $B_X$.
\par But the weak maximum principle holds.
\bT\label{T:mp} If $X$ is a compact set in ${\mathbb C}^n$ and $u\in
PSH(X)$, then $u(z)\le\sup_{w\in B_X} u(w)$ for all $z\in X$.
\eT
\begin{proof} By \cite[Thm. 4.2]{P2} for all points $z\in X$ there is
$L\in\T_z(X)$ with $\supp L\sbs B_X$. Hence
$u(z)\le\mu_L(u)\le\sup_{w\in B_X} u(w)$.
\end{proof}
\par A compact set $X$ in ${\mathbb C}^n$ is called {\em $O$-regular}
if the set $O_X$ is closed. As the following theorem shows such
sets admit continuous solutions in the class of plurisubharmonic
function for the Dirichlet problem with continuous boundary data
on $B_X$.
\bT \label{LemmaC} If $X$ is
$O$-regular, then for every $\phi\in C(B_X)$ there is a
continuous plurisubharmonic function $u$ on $X$ equal to $\phi$
on $B_X$.
\eT
\par We start the proof with an approximation of our solution $u$
\bL\label{LemmaD}
Under the assumptions of Theorem \ref{LemmaC}, for every
$\epsilon
>0$ there is a continuous plurisubharmonic function $u$ on $X$
such that
\[
\phi-\epsilon \leq u\leq \phi\qquad \text{ on } \quad B_X
\]
and $u\geq \min_{z\in B_X}\phi(z)$ on $X$.
\eL
\begin{proof} For $w\in B_X$ let $u_w$ be a continuous plurisubharmonic peak
function at $w$ on $X$. Since $\phi(z)-\phi(w)>-\eps/2$ on a
neighborhood $U$ of $w$, there is a constant $c_w>0$ such that the
function
\[
v_w(z) =\phi(w)+c_wu_w(z)-\epsilon/2\le \phi(z)
\]
on $B_X$. Then $v_w(w)=\phi(w)-\epsilon/2$.
\par Let $V_w=\{z\in B_X : v_w(z)>\phi(z)-\epsilon \}$.
Since $V_w$ is a non-empty relatively open set in $B_X$
containing $w$ there is a finite cover of $B_X$ by sets
$V_{w_j}$, $1\leq j\leq m$.  If $a=\min _{z\in B_X}\phi(z)$,
consider the function
\[
u(z)=\max\{a,v_{w_1}(z),\dots,v_{w_m}(z)\}, \qquad z\in X.
\]
Then $u$ is a continuous plurisubharmonic function on $X$ and
$a\leq u\leq \phi$ on $B_X$.  Moreover, if $z\in V_{w_j}$, then
$u(z)\geq v_{w_j}(z)\geq\phi(z)-\epsilon$.
\end{proof}
\par Now we turn to the proof of Theorem \ref{LemmaC}
\begin{proof} By taking $\epsilon=1/2$
in Lemma \ref{LemmaD} we can find a continuous plurisubharmonic
function $u_1$ on $X$ such that $\phi-1/2\leq u_1\leq \phi$ on
$B_X$. Let $\phi_1=\phi-u_1$ be a function on $B_X$.  Then $0\leq
\phi_1\leq 1/2$ on $B_X$. Again by Lemma \ref{LemmaD} we find a
continuous plurisubharmonic function $u_2$ on $X$ such that
$u_2\geq 0$ on $X$ and $\phi_1-1/4\leq u_2\leq \phi_1$ on $B_X$.
By the maximum principle $u_2\leq 1/2$ on $X$. Let
$\phi_2=\phi_1-u_2$.  Then $0\leq \phi_2\leq 1/4$.
\par Suppose that for $j=2,3,\dots,k$ the functions $u_j$ and
$\phi_j$ have been chosen  such that:

\be
\item $u_j$ are plurisubharmonic and continuous  and $0\leq u_j\leq 1/2^{j-1}$
on $X$.
\item $\phi_j$ are continuous, $0\leq \phi_j\leq 1/2^j$,
$\phi_{j-1}-1/2^j \leq u_j\leq \phi_{j-1}$, and
$\phi_j=\phi_{j-1}-u_j$   on $B_X$,
\ee
By Lemma \ref{LemmaD} we find a continuous plurisubharmonic
function $u_{k+1}$ on $X$ such that $u_{k+1}\geq 0$ on $X$ and
\[
\phi_{k}-1/2^{k+1}\leq u_{k+1}\leq\phi_k \qquad \text{ on } \quad B_X
\]
By the maximum principle $u_{k+1}\leq 1/2^k$ on $X$.  Let
$\phi_{k+1}=\phi_k-u_{k+1}$.  Then $0\leq \phi_{k+1}\leq
1/2^{k+1}$ on $B_X$.   Thus the functions with the properties (1)
and (2) are defined for all natural $j$.

Now we let
\[
u(z)=\sum_{j=1}^\infty u_j(z), \qquad z\in X.
\]
Then $u$ is continuous and plurisubharmonic on $X$. Moreover, for
every $k$ we have
\[
\phi-u=\phi-\sum_{j=1}^\infty u_j=\phi_k-\sum_{j=k+1}^\infty u_j
\]
on $B_X$. Since the right hand side tends to $0$ as $k\to \infty$
we have proved that $u=\phi$ on $B_X$.
\end{proof}
\section{Poisson sets}\label{S:ps}
\par A continuous real valued function $h$ on $X$ is said to be {\it
pluriharmonic} if $h(z)=\mu(h)$ for all $\mu\in \J_z(X)$. A
compact set is called a {\it Poisson set} if for every $\phi\in
C(B_X)$ there is a pluriharmonic function $h$ on $X$ such that
$h=\phi$ on $B_X$.
\bT A compact set $X\subset {\mathbb C}^n$ is a Poisson set if and
only if for every $z\in X$ the set $\J_z^b(X)$ contains only
one measure.  This measure will be denoted by $P_z$.
\eT

\begin{proof} Suppose that $X$ is Poisson, $z\in X$ and
$\mu_1,\mu_2\in \J_z(X)$. Let $\phi\in C(B_X)$ and let $u$ be a
pluriharmonic function on $X$ such that $u=\phi$ on $B_X$.
Then $u(z)=\mu_1(\phi)=\mu_2(\phi)$ and we conclude that
$\mu_1=\mu_2$.
\par Now suppose that for every $z\in X$ the set $\J_z^b(X)=\{P_z\}$.
Let $\phi\in C(B_X)$ and $h(z)=P_z(\phi)$.  Let us show that $h$
is continuous.  Suppose that the sequence $z_j$ in $X$ converges
to $z$. Switching if necessary to a subsequence, we may assume
that the measures $P_{z_j}$ converge weak-$*$ to some measure
$\mu\in \J^b_z(X)$. By assumption, $\mu=P_z$ and $h(z)=\lim_{j\to
\infty} P_{z_j}(h)=\lim_{j\to \infty} h(z_j)$.
\par Let us now prove that $h$ is plurisubharmonic.  By Theorem \ref{LemmaC} we
can find a continuous plurisuperharmonic function $v$ on $X$ equal
to $\phi$ on $B_X$.   The function
\[
E_Xv(z)=\inf \{\mu(v) : \mu\in \J_z(X)\}
\]
is lower semicontinuous, and $\mu(E_Xv)\geq E_Xv(z)$ for all
$\mu\in \J_z(X)$.
\par By Theorem \ref{T:lsp} for every $\mu\in \J_z(X)$ there is
$\nu\in \J_z^b(X)$ such that for every continuous plurisubharmonic
function $u$ on $X$ we have $\mu(u)\leq \nu(u)$. By assumption we
have $\J_z^b(X)=\{P_z\}$.  Hence $\mu(u)\leq P_z(u)$ for every
$\mu \in \J_z(X)$. Since $v$ is plurisuperharmonic,  $\mu(v)\geq
P_z(v)$. Hence $E_Xv(z)=P_z(v)=h(z)$. Since $h$ is continuous and
$h(z)\leq \mu(h)$ for all $\mu\in \J_z(X)$, the function $h$ is
plurisubharmonic.
\par Now we prove that $h$ is pluriharmonic.  If $\mu\in \J_z(X)$,
$z\in X$, then
\[
h(z)\leq \mu(h)\leq P_z(h)=P_z(\phi)=h(z)
\]
Hence $h(z)=\mu(h)$ and $h$ is pluriharmonic
\end{proof}
\section{Maximal solutions of the Dirichlet problem}\label{S:msdp}
\par In this section we prove the existence of maximal solutions for
the Dirichlet problem on a compact set $X\sbs{\mathbb C}^n$.
The Dirichlet problem for plurisubharmonic maximal functions is
to find a maximal function on $X$ equal to the given function
$\phi$ on $B_X$. Even if $\phi\in C(B_X)$ an upper
semicontinuous solution does not need to exist. For example, if
$X_1=\{(z_1,0):\,z_1\in\ovr{\mathbb D}\}$,
$X_2=\{(0,z_2):\,z_2\in\ovr{\mathbb D}\}$ and $X=X_1\cup
X_2\sbs{\mathbb C}^2$, then $B_X$ is the union of the sets
$B_{X_1}=\{(z_1,0):\,z_1\in{\mathbb T}\}$ and $B_{X_2}=
\{(0,z_2):\,z_2\in{\mathbb T}\}$. If $\phi\equiv1$ on $B_{X_1}$
and 0 on $B_{X_2}$, then a maximal solution of the Dirichlet
problem with the boundary values $\phi$ should be equal 1 on
$X_1$ and to $0$ on $X_2$. So at the origin we have to choose
between 0 and 1. If choose 1 we loose the subaveraging property
and if we choose 0 we loose upper semicontinuity.
\par It seems more reasonable to save the subaveraging property and
introduce the cone of {\em weakly plurisubharmonic} functions
on $X$ denoted by $PSH^w(X)$. A function belongs to $PSH^w(X)$
if it is the upper envelope of a family of continuous
plurisubharmonic functions on $X$. Such functions are lower
semicontinuous and satisfy the subaveraging inequality for
every $z\in X$ and every $\mu\in\J_z(X)$. In particular, the
envelopes of continuous functions are weakly plurisubharmonic.
\bL\label{L:dwpsh} A function $u\in PSH^w(X)$ if and only if
$u$ is the limit of an increasing sequence of continuous
plurisubharmonic functions defined on neighborhoods of $X$.
\eL
\begin{proof} We will prove only the necessity of this
condition because the sufficiency is trivial. Since $u$ is
lower semicontinuous, it is the limit of an increasing sequence
of continuous functions $\phi_j$ on $X$. For every $j$ and
$z\in X$ find a function $v_{jz}\in PSH^c(X)$ such that
$v_{jz}\le u$ and $v_{jz}\ge\phi_j-2/j$ on a non-empty
neighborhood $U_{jz}$ of $z$. Choose finitely many points
$z_{j1},\dots, z_{jk_j}$ such that the sets $U_{jz_{jk}}$ cover
$X$ and let $v_j=\max\{v_{jz_{jk}},\,1\le k\le k_j\}$. The
functions $v_j\in PSH^c(X)$, $v_j\le u$ and $v_j\ge\phi_j-2/j$.
\par By Corollary 3.1 from \cite{P2} each $v_j$ is the uniform limit of
an increasing sequence of continuous plurisubharmonic functions
defined on neighborhoods of $X$. So for each $j$  we can find a
continuous plurisubharmonic function $u_j$ defined on an open
neighborhood $Y_j$ of $X$ such that $u_j\le u$ and
$u_j\ge\phi_j-1/j$ on $X$. Clearly, the sequence $\tilde
u_j=\max\{u_1,\dots,u_j\}$ satisfies all requirements.
\end{proof}
\par Let $Z$ be a closed set in $X$ containing $B_X$. A  weakly
plurisubharmonic function $u$ on $X$ is {\em maximal} on $X\sm Z$
if for every relatively open $V\subset X\sm Z$ and any function
$v\in PSH^c(\ovr V)$ the inequality $v\leq u$ on the relative
boundary $\bd_XV$ of $V$ in $X$ implies that $v\leq u$ on
$\overline V$. If $Z=B_X$ then we say that $u$ is maximal.
\par It should be noted that our definition of maximal functions on
compact sets differs from the definition of maximal functions on
open sets. For example, if $X=\ovr{\mathbb D}^2$ then the
functions $u(z)\equiv0$ and $v(z)=|z_1|^2-1$ coincide on
$B_X=\mathbb T^2$ and both are maximal on $\mathbb D^2$ in the
sense that their Monge--Amp\`ere mass is equal to 0. However, the
function $v$ is not maximal on $X$ according to our definition.
\par The weakly plurisubharmonic maximal functions have the following
description.
\bT\label{T:cmf} Let $X\sbs{\mathbb C}^n$ be a compact set, let $Z$ be
a closed set in $X$ containing $B_X$ and let $u\in PSH^w(X)$.
Then $u$ is maximal on $X\sm Z$ if and only if for any $z_0\in
X$ there is a measure $\mu\in\J_{z_0}(X)$ such that
$\supp\mu\sbs Z$ and $\mu(u)=u(z_0)$.
\eT
\begin{proof} Suppose that $u$ is maximal on $X\sm Z$. Let $\phi$ be
a continuous function on $X$ equal to $0$ on $Z$ and greater
than 0 on $X\sm Z$. We define $v=E_X(u+\phi)$. Since $u+\phi\le
u$ on $Z$, $v\le u$. On the other hand, for every $z\in X$
\[v(z)=\inf\{\mu(u+\phi) : \mu\in \J_z(X)\}\ge E_Xu(z)=u(z).\]
Hence $v=u$. Edwards' Theorem now gives that there exists a
measure $\mu$ in $\J_{z_0}(X)$ such that $\mu(u+\phi)=u(z_0)$.
If $(X\sm Z)\cap\supp\mu$ is non-empty, then
$\mu(u+\phi)>\mu(u)\ge u(z_0)$. Thus $\supp\mu\sbs Z$ and
$u(z_0)=\mu(u+\phi)=\mu(u)$.
\par Now suppose that for any $z_0\in X$ there is
a measure $\mu\in\J_{z_0}(X)$ such that $\supp\mu\sbs Z$ and
$\mu(u)=u(z_0)$.  To show that $u$ is maximal on $X\sm Z$, we
consider a relatively open subset $V$ of $X\sm Z$ and a continuous
plurisubharmonic function $v$ on $\overline V$ such that $v\leq u$
on $\bd_X V$. Then we assume that $z_0\in V$ and take a totally
perfect sequence $L=\{f_j\}\in \T_{z_0}(X)$ with $\mu_L=\mu$.
\par   Let $W$ be an open set in ${\mathbb C}^n$
such that $V=X\cap W$ and $\ovr V=X\cap\ovr W$. Let $D_j$ be the
connected component of $f_j^{-1}(W)$ containing the origin.  We
let $a_j:{\mathbb D}\to D_j$ be a holomorphic universal covering
with $a_j(0)=0$ and consider the sequence $M=\{g_j\}$,
$g_j=f_j\circ a_j$. Clearly $\cl M\subset\ovr V$ and $z_M=z_0$.
\par Switching to a subsequence if necessary we may assume that the
measures $\mu_{g_j}$ converge weak-$*$ to a measure $\mu_M$. We
claim that $\supp\mu_M\subset\partial_X V$.
\par To prove it we take a relatively open set $U\sbs V$ such that
$\ovr U\sbs V\sm \bd_XV$. We also take an open set $Y$ in
${\mathbb C}^n$ such that $U=Y\cap X$ and $Y\subset \subset W$.
Let $E_j$ be the set of all $\zeta\in\mathbb T$ where $a_j$ has a
radial limit.  We will denote this limit by $a_j^*(\zeta)$.  Let
$F_j\subset E_j$ be the set where $g_j$ has a radial limit and
this limit belongs to $\overline Y$. Our goal is to prove that
$\lambda(F_j)\to 0$. Since $\mu_{g_j}=\lambda(F_j)$ it will show
that $\mu_{g_j}(\overline Y)\to 0$ and, consequently,
$\mu_M(\overline Y)=0$. Thus $\supp\mu_M\subset\partial_X V$.
\par Our first observation is that if $\zeta\in F_j$ and
$a_j^*(\zeta)=\xi\in {\mathbb D}$, then $f_j(\xi)\in \overline
Y\subset\subset W$ and, therefore, $\xi\in D_j$. But $a_j$ is a
universal covering and, therefore, $a_j^*(\zeta)\in \partial D_j$.
and this contradiction shows that $a_j^*(\zeta)\in {\mathbb T}$.
Hence the set $G_j=a_j^*(F_j)$ lies in $\mathbb T$.
\par If $\xi\in G_j$ then we can find $\zeta\in F_j$ such that
the path $a_j(r\zeta)$, $0\leq r<1$, lies in ${\mathbb D}$ and
ends at $a_j^*(\zeta)=\xi\in\mathbb T$.  Since $g_j$ has a
limit along the path $r\zeta$, $0\le r<1$, $f_j$ has a limit
along the path $a_j(r\zeta)$, $0\leq r<1$. By Lindel\"of's
theorem, $f_j$ has the radial limit at $\xi$ and $\lim_{r\to
1^-}f_j(r\xi)=\lim_{r\to 1^-}g_j(r\zeta)$. So, we have shown
that $f_j$ has a radial limit at all points of $G_j$ and
$f_j^*(G_j)\subset \overline Y$. Since $\overline Y\cap
Z=\emptyset$ and $\supp \mu\subset Z$, $\mu_{f_j}(\overline
Y)\to 0$.  Hence $\lambda(G_j)\to 0$.
\par Let us show that also $\lambda(F_j)\to 0$. For this we take
an open set $H_j\sbs\mathbb T$ such that $G_j\subset H_j$ and
$\lambda(H_j)\to 0$.  Let $h_j(\xi)=\omega(\xi,H_j,{\mathbb
D})$ and $\tilde h_j(\zeta)=h_j(a_j(\zeta))$ for $\zeta\in
{\mathbb D}$. Then for every $\zeta_0\in F_j$ we have
\[
\lim_{r\to 1^-}\tilde h_j(r\zeta_0) =\lim_{r\to
1^-}h_j(a_j(r\zeta_0))=1,
\]
because $a_j^*(\zeta_0)\in G_j$.  Thus,
$\omega(\zeta,F_j,{\mathbb D})\leq\tilde h_j(\zeta)$.  Since
$a_j(0)=0$,
\[
\lambda(F_j)=\omega(0,F_j,{\mathbb D})\leq \tilde
h_j(0)=h_j(0)=\lambda(H_j).
\]
Thus $\lambda(F_j)\to 0$ and $\supp \mu_M\subset\partial_XV$.

\par By Lemma \ref {L:dwpsh} for any $\eps>0$ we can find a
continuous plurisubharmonic function $w$ defined on a
neighborhood of $X$ such that $w\le u$ on $X$ and
$\mu_M(u)\le\mu_M(w)+\eps$. So we may consider functions
$w_j=w\circ f_j$ when $j$ is large.
\par By Theorem \ref{T:lmp}
\[
\mu_{g_j}(w)=\dfrac 1{2\pi}\int_0^{2\pi}w_j(a_j(e^{i\theta}))\,
d\theta \leq \dfrac 1{2\pi}\int_0^{2\pi}w_j(e^{i\theta})\, d\theta
= \mu_{f_j}(w)
\]
and we see that $\mu_M(w)\leq \mu_L(w)$. Since $\mu_M\in
\J_{z_0}(\ovr V)$, we have $v(z_0)\leq \mu_M(v)\leq \mu_M(u)$.
Therefore,
\[v(z_0)\le\mu_M(u)\le\mu_M(w)+\eps\le\mu_L(w)+\eps=u(z_0)+\eps.
\]
Since $\eps>0$ is arbitrary we see that $v(z_0)\le u(z_0)$.
Hence $u$ is maximal.
\end{proof}
\par So in the example at the beginning of this section the
function $u$ equal to 0 on $X_2$ and 1 on $X_1\sm\{0\}$ is
maximal and weakly plurisubharmonic.

\par The theorem below establishes the existence of maximal
solutions of the Dirichlet problem which satisfy the
subaveraging inequality but are lower semicontinuous.
\bT\label{TheoremE}
Let $X$ be an $O$-regular set and let $\phi$ be a continuous
function on $B_X$. Define $u$ by
\[
u=\sup\{v: v\in PSH^c(X), v|_{B_X}\leq \phi\}.
\]
Then $u$ has the following properties:
\be
\item $u$ is lower semicontinuous and weakly plurisubharmonic on $X$;
\item for every $z\in X$
\[u(z)=\inf\{\mu(\phi) : \mu\in \J_z^b(X)\};\]
\item for every $z\in X$ there is $\mu\in \J_z^b(X)$ such that
$u(z)=\mu(u)$;
\item $u=\phi$ on $B_X$ and $\lim_{z\to z_0}u(z)=\phi(z_0)$ for every
$z_0\in B_X$;
\item $u$ is maximal on $X$;
\item $u$ is the only weakly plurisubharmonic maximal on $X$ function
equal to $\phi$ on $B_X$.
\ee
\eT
\begin{proof} (1): Since  $u$ is a supremum of a family of continuous
functions, $u$  is lower semicontinuous and weakly
plurisubharmonic on $X$.
\par (2): Let $u_1(z)=\inf\{\mu(\phi) : \mu\in \J_z^b(X)\}$. It is
obvious that $u(z)\le u_1(z)$ because every $v\in PSH^c(X)$
such that $v|_{B_X}\le\phi$ satisfies this inequality. By
Theorem \ref{LemmaC}  there is a continuous plurisuperharmonic
function $u_2$ on $X$ equal $\phi$ on $B_X$. Let
\[u_3=\sup\{v: v\in PSH^c(X), v\leq u_2\}.\]
Clearly, $u_3\le u$. By Edwards' Theorem $u_3(z)=\inf\{\mu(u_2)
: \mu\in \J_z(X)\}$. By Theorem \ref{T:lsp} for every $\mu\in
\J_z(X)$ there is a measure $\nu\in \J_z^b(X)$ such that
$\mu(u_2)\geq \nu(u_2)$. Hence
\[
u_3(z)=\inf\{\mu(u_2) : \mu\in \J_z^b(x)\}=\inf\{\mu(\phi) :
\mu\in \J_z^b(x)\}=u_1(z).
\]
Thus $u\equiv u_1$.
\par (3): Since the set $\J_z^b(X)$ is weak-$*$ compact, by (2) there is
$\mu_z\in \J_z^b(X)$ such that $u(z)=\mu_z(\phi)$.
\par (4): Suppose $z_0\in B_X$ and take a sequence $\{z_j\}$ in $X$
converging to $z_0$.   Since $\J_{z_0}(X)=\{\delta_{z_0}\}$, by
(3) $u(z_0)=\phi(z_0)$.  Let us take $\mu_j\in \J_{z_j}^b(X)$ such
that $\mu_{z_j}(v)=u(z_j)$.  The sequence $\{\mu_j\}$ must
converge weak-$*$ to $\delta_{z_0}$, otherwise $\J_{z_0}(X)$ will
contain a measure different from $\delta_{z_0}$.  Hence $\lim
u(z_j)=\phi(z_0)=u(z_0)$.
\par (5): Follows immediately from Theorem \ref{T:cmf}.
\par (6): Let $v$ be a weakly plurisubharmonic maximal function
on $X$ equal to $\phi$ on $B_X$. Since by (3) for every $z\in
X$ there is $\mu\in\J^b_z(X)$ such that $u(z)=\mu(u)$,
$v(z)\le\mu(\phi)=u(z)$. Thus $v\le u$.
\par On the other hand since $v$ is maximal on $X$, $v$ is greater
or equal than any continuous plurisubharmonic function which does
not exceed $\phi$ on $B_X$. But $u$ is the upper envelope of the
latter functions and this implies that $u\le v$. Hence $u=v$.
\end{proof}
\par We will call the function $u$ from the theorem above {\em the
maximal solution of the Dirichlet problem with a boundary value
$\phi$}.
\section{Harmonicity of maximal functions}\label{S:homf}
\par Given a weak-$*$ converging sequence $L=\{f_j\}$ and a point
$z\in\cl L$ we denote by $\J^s_z(L)$ all measures $\mu$ for
which there are a subsequence $\{j_k\}$ and a sequence of
holomorphic mappings $g_k:\,{\mathbb D}\to{\mathbb D}$ such
that the measures $\mu_{h_k}$, $h_k=f_{j_k}\circ g_k$, converge
weak-$*$ to $\mu$ and points $h_k(0)$ converge to $z$. Let
$\J_z(L)$ be the convex envelope of $\J^s_z(L)$. Clearly,
$\J_z(L)\sbs\J_z(\cl L)$. A totally perfect sequence $L$ is
called {\it ample} if $\J_z(L)=\J_z(\cl L)$ for all $z\in\cl
L$.
\par Let us denote by $PSH^c(L)$ the set of continuous
function $u$ on $\cl L$ such that $u(z)\le\mu(u)$ for every
$z\in\cl L$ and every $\mu\in\J_z(L)$. A sequence $L$ is ample if
and only if $PSH^c(L)=PSH^c(\cl L)$. Indeed, the necessity is
clear and the sufficiency follows from the Hahn--Banach and
Edwards' theorems. Indeed, if there is $\mu\in\J_z(\cl
L)\sm\J_z(L)$ then there is a continuous function $\phi$ on $\cl
L$ such that
\[\mu(\phi)<\inf_{\nu\in\J_z(L)}\nu(\phi).\] By Edwards' theorem
the infimum above is equal to the value at $z$ of the envelope
of $\phi$ in $PSH^c(L)$ which, in turn, is equal to $E_{\cl
L}\phi(z)$. But the latter value does not exceed $\mu(\phi)$
and we got a contradiction.
\par  Let $L=\{f_j\}$ be a perfect sequence. A continuous plurisubharmonic function
$u$ on $X=\cl L$ is said to be {\em harmonic on} $L$ if
$\mu(u)=u(z)$ for every $z\in\cl L$ and every $\mu\in\J_z(L)$. A
continuous plurisubharmonic function $u$ on a compact set $X$ is
said to be {\em harmonic} if $\mu(u)=u(z)$ for every $z\in X$ and
every $\mu\in\J_z(X)$.
\par The theorem below relates harmonicity and minimizing
sequences. But before we have to prove a technical lemma.
\bL\label{L:sp} Let $L=\{f_j\}$ be a totally perfect sequence in $\mathbb
C^n$ and let $\{A_j\}$ be a sequence of closed sets in $\mathbb
D\sm\{0\}$ such that $\om_j=\om(0,A_j,\mathbb D)\to0$ as
$j\to\infty$. Let $p_j$ be a holomorphic covering mapping of
$D_j=\mathbb D\sm A_j$ such that $p_j(0)=0$. Then the sequence
$L'=\{f_j\circ p_j\}$ is totally perfect, $z_{L'}=z_L$, $\cl
L'=\cl L$, and $\mu_{L'}=\mu_L$.
\eL
\begin{proof} Let $g_j=f_j\circ p_j$. We start with proving that $\mu_{L'}=\mu_L$. If
$\phi\in C(\mathbb C^n)$ then $\mu_{f_j}(\phi)=u_j(0)$, where
$u_j$ is a harmonic function on $\mathbb D$ equal to
$\phi_j=\phi\circ f_j$ on $\mathbb T$, while
$\mu_{g_j}(\phi)=v_j(0)$, where $v_j$ is a harmonic function on
$\mathbb D\sm A_j$ equal to $\phi_j$ on $\mathbb T\cup\bd A_j$.
Hence
\[|\mu_{f_j}(\phi)-\mu_{g_j}(\phi)|=|u_j(0)-v_j(0)|\le\|\phi_j\|\om_j.\]
Since the sequence $L$ is uniformly bounded we see that
$\mu_{L'}=\mu_L$.
\par Clearly, $z_{L'}=z_L$ and $\cl L'\sbs\cl L$. If $z\in\cl L$,
$r>0$ and $\eps_j=\om(0,f_j^{-1}(B(z,r)),\mathbb D)$, then
$\eps_j\ge a>0$. Now the same reasoning as above shows that
$\om(0,g_j^{-1}(B(z,r)),\mathbb D)\ge a-\om_j$ and this tells
us that $z\in\cl L'$ and $L'$ is totally perfect.
\end{proof}
\par Now we can prove the theorem.
\bT\label{T:hos} Let $X\sbs{\mathbb C}^n$ be a compact set and $u\in PSH^c(X)$
be a plurisubharmonic function. If $L=\{f_j\}\in\T_{z_0}(X)$
and $u(z_0)=\mu_L(u)$, then there is a sequence $L'=\{g_j\}$ in
$\T_{z_0}(X)$ such that $\cl L'=\cl L$, $\mu_{L'}=\mu_L$ and
$u$ is harmonic on $L'$.
\par Moreover, there is a sequence $\{H_j\}$ of harmonic functions
on $\mathbb D$ such that
\[\lim_{m\to\infty}H_{j_m}(\zeta_m)=u(z)\]
for any subsequence $j_m$ and any sequence
$\{\zeta_m\}\sbs\mathbb D$ such that the points
$g_{j_m}(\zeta_m)$ converge to $z\in \cl L'$
\eT
\begin{proof}
\par By Lemma 3.1 from \cite{P2} there is an increasing sequence of
continuous plurisubharmonic functions $\{u_k\}$ on neighborhoods
$V_k$ of $X$ converging pointwise to $u$ on $X$. It follows that
the functions $u_k$ converge to $u$ uniformly on $X$. So if $\phi$
is a continuous extension of $u$ to ${\mathbb C}^n$, then we can
assume that $|u_k-\phi|<1/k$ on $V_k$.
\par For every $k$ we consider the functions $u_{jk}=u_k\circ f_j$
which are defined when $j$ is large.  Let $F_{jk}$ be the
harmonic functions on $\mathbb D$ equal to $u_{jk}$ on
${\mathbb T}$. Then $F_{jk}\geq u_{jk}$ and
$F_{jk}(0)=\mu_{f_j}(u_k)$.
\par Let $\dl_k=|u(z_0)-\mu_L(u_k)|$. Since
$$\lim_{k\to\infty}\mu_L(u_k)=\mu_L(u)=u(z_0),$$ we see that the
sequence $\{\dl_k\}$ converges to 0. Since
$\lim_{j\to\infty}\mu_{f_j}(u_k)=\mu_L(u_k)$ for each $k$ we can
find $j_k$ such that $|\mu_{f_{j_k}}(u_k)-\mu_L(u_k)|<1/k$. Hence
\begin{equation}\begin{aligned}
0&\le F_{j_kk}(0)-u_k(z_0)\notag\\&=
\mu_{f_{j_k}}(u_k)-\mu_L(u_k)+\mu_L(u_k)-u(z_0)+u(z_0)-u_k(z_0)
\notag\\&\le\frac2k+\dl_k=\gm_k.\notag
\end{aligned}\end{equation}
\par Let $v_k(\zeta)=u_k(f_{j_k}(\zeta))-F_{j_kk}(\zeta)$ and let
$A_k=\{\zeta\in{\mathbb D}:v_k(\zeta)\le-\sqrt{\gm_k}\}$. The
functions $v_k$ are negative and subharmonic on ${\mathbb D}$
and $v_k(0)\ge-\gm_k$. By Two Constants Theorem
$v_k(0)\le-\sqrt{\gm_k}\om(0,A_k,{\mathbb D})$, i.e.,
$\om(0,A_k,{\mathbb D})\le\sqrt{\gm_k}$.
\par Let $D_k={\mathbb D}\sm A_k$ and let $p_k$ be a holomorphic
covering mapping of $D_k$ by ${\mathbb D}$ such that
$p_k(0)=0$. We let $g_k=f_{j_k}\circ p_k$. If $w_k=u_k\circ
g_k$ and $F_k=F_{j_kk}\circ p_k$, then
\[0\ge w_k-F_k\ge-\sqrt{\gm_k}\] on ${\mathbb D}$. Hence, if $H_k$ is
the harmonic functions on $\mathbb D$ equal to $w_k$ on
${\mathbb T}$, then
\[0\ge H_k-F_k\ge w_k-F_k\ge-\sqrt{\gm_k}\] on ${\mathbb D}$.
Therefore,
\begin{equation}\label{e:hos1}
0\ge w_k-H_k=w_k-F_k+F_k-H_k\ge-2\sqrt{\gm_k}.
\end{equation}
\par We let $L'=\{g_k\}$. By Lemma \ref{L:sp}, $L'\in\T_{z_0}(X)$ is
totally perfect, $\mu_{L'}=\mu_L$ and $\cl L=\cl L'$. If $z\in\cl
L$ and $\nu\in\J^s_z(L)$, then there are a subsequence $\{k_m\}$
and a sequence of holomorphic mappings $q_m:\,{\mathbb
D}\to{\mathbb D}$ such that $\nu=\mu_M$ and $z=z_M$, where
$M=\{h_m\}=\{g_{k_m}\circ q_m\}$. Let $\zeta_m=q_m(0)$. Consider
the mappings
\[e_m(\zeta)=
g_{k_m}\left(\frac{\zeta+\zeta_m}{1+\ovr\zeta_m\zeta}\right).\]
By the Littlewood subordination principle
$\mu_{e_m}(w_{k_m})\ge\mu_{h_m}(u_{k_m})$. Since
$\mu_{e_m}(u_{k_m})=H_{k_m}(\zeta_m)$, by (\ref{e:hos1})
\begin{equation}\begin{aligned}
\mu_{h_m}(u_{k_m})&\le\mu_{e_m}(u_{k_m})=H_{k_m}(\zeta_m)\le
w_{k_m}(\zeta_m)+2\sqrt{\gm_k}\notag\\&=u_{k_m}(h_m(0))+2\sqrt{\gm_k}\le
\phi(h_m(0))+2\sqrt{\gm_k}+\frac1{k_m}.\notag
\end{aligned}\end{equation}
Note that $|\mu_{h_m}(u_{k_m})-\mu_{h_m}(\phi)|\le1/k_m$. So
taking the limits as $m\to\infty$ in the inequality above we
get
\[\mu_M(u)=\lim_{m\to\infty}\mu_{h_m}(u_{k_m})
\le\liminf_{m\to\infty}H_{k_m}(\zeta_m)\le
\limsup_{m\to\infty}H_{k_m}(\zeta_m)\le u(z_M).\] But $\mu_M\in\J_{z_M}(X)$
and, therefore,
\[\mu_M(u)=\lim_{m\to\infty}\mu_{h_m}(u_{k_m})=
\lim_{m\to\infty}H_{k_m}(\zeta_m)=u(z_M).\] Since this equality
holds for all $\mu\in\J^s_z(L')$, $z=z_M$, it holds for all
measures in $\J_z(L')$. Hence $u$ is harmonic on $L'$. Note
that the last equality also proves the second claim of the
theorem.
\end{proof}
\par The following corollary follows immediately from Theorems
\ref{TheoremE} and \ref{T:hos} and the definition of ample
sequences.
\bC Let $X$ be a $O$-regular compact set and $u\in C(X)$ be a maximal
plurisubharmonic function. Then for every $z_0\in X$ there is a
sequence $L\in\T_{z_0}(X)$ such that $\mu_L\in\J^b_{z_0}(X)$
and $u$ is harmonic on $L$. Moreover, if $L$ is ample then $u$
is harmonic on $\cl L$.
\eC
\par By Theorem \ref{T:cmf} the values of maximal functions are
determined by measures $\mu$ in $\J_z(X)$ minimizing the
functional $\mu(\phi)$. For any such measure $\mu$ there is a
totally perfect sequence $L\in\T_z(X)$ with $\mu_L=\mu$. As the
following theorem shows the function $u$ stays maximal on $\cl
L$.
\bT\label{T:hamf} Let $X$ be a compact set in $\mathbb C^n$, let
$Z$ be a closed set in $X$ containing $B_X$ and let $u\in
PSH^c(X)$ be a function maximal on $X\sm Z$. If $z_0\in X\sm
Z$, $\mu\in\J_{z_0}(X)$, $\supp\mu\sbs Z$, $\mu(u)=u(z_0)$ and
$L=\{f_j\}$ is a totally perfect sequence such that $z_L=z_0$
and $\mu_L=\mu$, then $u$ is maximal on $\cl L\sm Z$ and
$B_{\cl L}\sbs Z$.
\eT
\begin{proof} Let $Y=\cl L$. By Lemma 2.2 in \cite{P2} for every
$z\in Y$ there is a totally perfect sequence $M$ such that
$z_M=z$, $\cl M\subset \cl L$, and $\supp \mu_M\subset \supp
\mu_L\subset Z$. Hence, if $z\not\in\supp\mu_L$ then
$\J_z(Y)\neq\{\dl_z\}$ and, consequently, $z\not\in O_Y$. Thus
$B_Y\sbs Z$.
\par To prove that $u$ is maximal on $Y\sm Z$ we replace $L$ with
$L'=\{g_j\}$ from Theorem \ref{T:hos}. It does not change
neither clusters nor supports. If $z\in Y\sm Z$ and
$\{\zeta_j\}$ is a sequence in $\mathbb D$ such that
$g_j(\zeta_j)\to z$ as $j\to\infty$, then the sequence of
mappings
\[h_j(\zeta)=
g_j\left(\frac{\zeta+\zeta_j}{1+\ovr\zeta_j\zeta}\right)\]
contains a weak-$*$ converging subsequence $M$ such that
$z_M=z$, $\mu_M\in\J_z^s(L')$ and $\supp\mu_M\sbs Z$. By
Theorem \ref{T:hos} $\mu_M(u)=u(z)$ and by  Theorem \ref{T:cmf}
it means that $u$ is maximal on $Y\sm Z$.
\end{proof}
\par The example at the beginning of the previous section and continued after
the proof of Theorem \ref{T:cmf} explains why we had to require
the continuity of $u$ in the theorem above. If $z_0=(1/2,0)$ we
can take as $L$ a constant sequence
\[L=\left\{\left(\frac{\zeta+1/2}{1+\zeta/2},0\right)\right\}.\]
This sequence satisfies all requirements of Theorem
\ref{T:hamf} but $\cl L=X_1$ and $u$ is not maximal on $X_1$
because $u(z,0)=1$, $z\ne0$, and $u(0,0)=0$.
\par The following theorem summarizes  our results regarding maximal functions
and totally perfect sequences. The part (1) was proved in Theorems
\ref{T:cmf} and \ref{T:hamf}. The second part asserts the
existence of a {\it minimal } sequence satisfying conditions of
Theorem \ref{T:cmf}.
\bT Let $X$ be a compact set in $\mathbb C^n$, let
$Z$ be a closed set in $X$ containing $B_X$ and let $u\in
PSH^c(X)$ be a function maximal on $X\sm Z$. For every point
$z_0\in X\sm Z$ there is a totally perfect sequence $L$ in $X$
with the following properties:\be
\item $z_L=z_0$, $\mu_L(u)=u(z_0)$, $\supp\mu_L\sbs Z$ and $u$ is
    maximal on $\cl L$;
\item if $M\in\T(\cl L)$, $z_M=z_0$, $\supp\mu_M\sbs Z$ and
    $\mu_M(u)=u(z_0)$, then $\cl L=\cl M$.\ee
\eT
\begin{proof} We will prove the second part. Let $\rS$ be the class of totally
perfect sequences satisfying (1). We introduce a partial order on
this set: $L\succeq M$ if $\cl M\sbs\cl L$. Let $\A$ be a totally
ordered subset of $\rS$. Let us show that $\A$ has a minimal
element in $\rS$. For every $L\in\A$ let $U_L=\mathbb C^m\sm\cl
L$. The sets $U_L$ are open and by Lindel\"of's Theorem (see
\cite[Thm. 1.4.14]{DS}) there is a sequence of sets $U_{L_j}$,
$j=1,\dots$, such that $\cup_{L\in\A}U_L=\cup_{j=1}^\infty
U_{L_j}$.
\par If $z\in U_{L_j}$ then $z\not\in\cl L_j$ and
\[\liminf_{i\to\infty}\Om(z_0,V,L_i)=0\]
for any open set $V$ containing $z$. By Theorems 2.3 and 2.2 from
\cite{P2} there is a totally perfect sequence $L$ such that
$z_L=z_0$, $\mu_L(u)=u(z_0)$, $\supp\mu_L\sbs Z$ and
\[\cl L=\{z\in\mathbb C^n:\,\liminf_{i\to\infty}\Om(z_0,V,L_i)>0\}.\]
Hence, $\cl L$ lies in $\mathbb
C^n\sm\cup_{L'\in\A}U_{L'}=\cap_{L'\in\A}\cl L'$ and this means
that $L\preceq M$ for any $M\in\A$. Moreover, by Theorem
\ref{T:hamf} $u$ is maximal on $\cl L$. Thus $L\in\rS$ and is a
minimal element for $\A$. By Zorn's lemma the class $\rS$ has a
minimal element $L$.
\par If $M\in\T(\cl L)$, $z_M=z_0$, $\supp\mu_M\sbs Z$ and
$\mu_M(u)=u(z_0)$, then $M\preceq L$ and by Theorem \ref{T:hamf}
$M\in\rS$. But $L$ is minimal and, therefore, $\cl L=\cl M$.
\end{proof}


\begin{thebibliography}{999}
\bibitem[BK]{BK} E. Bedford, M. Kalka, {\em Foliations and
    complex Monge-Amp\`ere equations,}. Comm. Pure Appl. Math. {\bf 30}(1977),
    543--571.
\bibitem[C]{C} J. Conway, {\em A Course in Functional
    Analysis,} Springer Verlag, 1990
\bibitem[D]{D} R. Dujardin, {\em Wermer examples and
    currents,} arXiv:0904.4179
\bibitem[DL]{DL} J. Duval, N. Levenberg, {\em Large polynomial
    hulls with no
    analytic structure,} Complex analysis and geometry, (Trento,
    1995), 119--122.
\bibitem[DS]{DS} N. Dunford, J. T. Schwartz, {\em Linear
Operators, Part I,} Interscience Publishers, 1958
\bibitem[Ga]{Ga} J. B. Garnett, {\em Bounded Analytic Functions,} Academic
Press, 1981
\bibitem[Go]{Go} G. M. Goluzin, {\em Geometric theory of functions of a
complex variable,} Translations of Mathematical Monographs, {\bf
26,} American Mathematical Society, Providence, R.I. 1969
\bibitem[LS]{LS} F. L\'arusson, R. Sigurdsson, {\em
    Plurisubharmonic functions and analytic disks on
    manifolds,} J. reine angew. Math., {\bf 501}(1998), 1-39.
\bibitem[L]{L} J. E. Littlewood, {\em On inequalities in the
    theory of functions,} Collected papers of J.E.
    Littlewood, vol. II, The Clarendon Press, Oxford University Press,
    New York, 1982, 963--1004.
\bibitem[Ph]{Ph} R. R. Phelps,  {\em Lectures on Choquet's
    Theory,}
    Lect.
    Notes. Math., {\bf 1757}, Springer Verlag, 2001
\bibitem[P1]{P1} E.A. Poletsky, {\em Holomorphic currents,}
Indiana Univ. Math. J., {\bf 42}(1993), 85--144
\bibitem[P2]{P2} E.A. Poletsky,  {\em Analytic geometry on
    compacta
    in ${\mathbb C}^n$,} Math. Z., {\bf 222}(1995), 407--424
\end{thebibliography}
\end{document}